\documentclass[12pt]{article}
\usepackage{amsfonts}
\usepackage{amssymb}
\usepackage{amsmath}
\usepackage{eufrak}
\usepackage[dvips]{graphicx}
\usepackage{latexsym}
\usepackage{verbatim}
\newcommand{\C}{{\EuFrak C}}

\newcommand{\Z}{{\mathbb{Z}}}

\newcommand{\bi}{\begin{itemize}}
\newcommand{\ei}{\end{itemize}}

\DeclareMathOperator{\GL}{GL}

\newtheorem{theorem}{Theorem}[section]
\newtheorem{lemma}[theorem]{Lemma}
\newtheorem{fact}[theorem]{Fact}
\newtheorem{claim}{Claim}
\newtheorem{subclaim}{Subclaim}

\newtheorem{definition}[theorem]{Definition}
\newtheorem{remark}[theorem]{Remark}

\newtheorem{maintheorem}{Theorem}

\title{On $\omega$-categorical, generically stable groups and rings}
\author{Jan Dobrowolski and Krzysztof Krupi\'nski\footnote{Research supported by the Polish Government grant N N201 545938}}
\date{}

\begin{document}
\maketitle
\begin{abstract} We prove that every $\omega$-categorical, generically stable group is nilpotent-by-finite and that every $\omega$-categorical, generically stable ring is nilpotent-by-finite.
\end{abstract}
\footnotetext{2010 Mathematics Subject Classification: 03C45, 03C35, 20A15.}
\footnotetext{Key words and phrases: $\omega$-categorical group, $\omega$-categorical ring, generic stability}

\section{Introduction}
A general motivation
 is to understand the structure of $\omega$-categorical groups and rings satisfying various natural model-theoretic assumptions. 

There is a long history of results of this kind. The fundamental theorem of Baur, Cherlin and Macintyre, proved in \cite{1}, says that $\omega$-categorical, stable groups are nilpotent-by-finite. A long-standing conjecture states that they are even abelian-by-finite,
which is known to be true in the superstable case. As to the $\omega$-categorical, stable rings, they are nilpotent-by-finite \cite{13}, and it is conjectured that they are null-by-finite. As for groups, this conjecture is known to be true in the 
superstable case. 
%
%
There are many generalizations and variants of these results. For example, $\omega$-categorical groups with NSOP (the negation of the strict order property) are nilpotent-by-finite \cite{6}, and $\omega$-categorical rings with NSOP are nilpotent-by-finite \cite{12}, too. It is also known that $\omega$-categorical, supersimple groups are finite-by-abelian-by-finite \cite{EW}, and $\omega$-categorical, supersimple rings are finite-by-null-by-finite \cite{KrWa}. 

More recently, in \cite{5}, an 
analysis of $\omega$-categorical groups and rings in the NIP environment has been undertaken. It was proved there that $\omega$-categorical rings with NIP are nilpotent-by-finite, and it was conjectured that $\omega$-categorical groups with NIP are nilpotent-by-finite, too. The conjecture was confirmed, but under the additional assumption of fsg (finitely satisfiable generics). It turns out that $\omega$-categorical groups with NIP and fsg are generically stable in the sense of \cite[Definition 6.3]{7}. This notion fits very well with the recent trend in model theory of studying structures whose some `pieces' are similar to stable ones. So, one can ask about the structure of $\omega$-categorical, generically stable groups and rings, without the NIP assumption. The following theorem was proved in \cite{11}.

\begin{fact}\label{solvable}
Every $\omega$-categorical, generically stable group is solvable-by-finite. 
\end{fact}

In this paper, using the above theorem, we obtain the following two results (to be more precise, in the proof of Theorem \ref{Thm2}, we will use an appropriate variant of Fact \ref{solvable}). 

\begin{maintheorem}\label{Thm1}
Every $\omega$-categorical, generically stable group is nilpotent-by-finite. 
\end{maintheorem}

\begin{maintheorem}\label{Thm2}
Every $\omega$-categorical, generically stable ring is nilpotent-by-finite. 
\end{maintheorem}

Theorem \ref{Thm1} proves \cite[Conjecture 3.5]{11}, and Theorem \ref{Thm2} answers \cite[Question 3.6]{11} in the affirmative.



\section{Preliminaries}

Recall that a first order structure $M$ in a countable language is said to be $\omega$-categorical if, up to isomorphism, $Th(M)$ has at most one model of cardinality $\aleph_0$. By Ryll-Nardzewski's theorem, this is equivalent to the condition that for every natural number $n$ there are only finitely many $n$-types over $\emptyset$. Assume $M$ is $\omega$-categorical. If $M$ is countable or a monster model (i.e. a model which is $\kappa$-saturated and strongly $\kappa$-homogeneous for a big cardinal $\kappa$), two finite tuples have the same type over $\emptyset$ iff they lie in the same orbit under the action of the automorphism group of $M$, and hence for each natural number $n$ the automorphism group of $M$ has only finitely many orbits on $n$-tuples (which implies that $M$ is locally finite). Moreover, 
for any finite subset $A$ of such an $M$, a subset $D$ of $M$ is $A$-invariant iff $D$ is $A$-definable.

 Let $T$ be a first order theory. We work in a monster model $\C$ of $T$.

Let $p \in S(\C)$ be invariant over $A\subset \C$. We say that $(a_i)_{i \in \omega}$ is a Morley sequence in $p$ over $A$ if 
$a_i \models p |A a_{<i}$ for all $i$. Morley sequences in $p$ over $A$ are indiscernible over $A$ and they have the same order type over $A$.
If $\C' \succ \C$ is a bigger monster model, then the generalized defining scheme of $p$ determines a unique $A$-invariant extension $\widetilde{p}\in S(\C')$ of $p$ (by the generalized defining scheme of $p$ we mean a family of sets $\{ p_i^\varphi: i \in I_\varphi\}$ (with $\varphi(x,y)$ ranging over all formulas without parameters) of complete types over $A$ such that $\varphi(x,c) \in p$ iff $c \in \bigcup_{i \in I_\varphi} p_i^{\varphi}(\C)$). By a Morley sequence in $p$ we mean a Morley sequence in $\widetilde{p}$ over $\C$. Finally, $p^{(k)}$ (where $k \in \omega \cup\{ \omega \}$) denotes the type over $\C$ of a Morley sequence in $p$ of length $k$.

\begin{definition}
A global type $p \in S(\C)$ is said to be generically stable if, for some small $A$, it is $A$-invariant and for each formula $\varphi(x;y)$ there is a natural number $m$ such that for any Morley sequence $(a_i: i< \omega)$ in $p$ over $A$ and any $b$ from $\C$ either less than $m$ $a_i$'s satisfy $\varphi(b;y)$ or less than $m$ $a_i$'s satisfy $\neg \varphi(b;y)$. 
In this definition, as a witness set $A$ one can take any (small) set over which $p$ is invariant. We will say that $p$ is generically stable over $A$ to express that $p$ is invariant over $A$ and generically stable.
\end{definition}

Recall \cite[Proposition 2.1]{8}.

\begin{fact}\label{indiscernible}
If $p$ is generically stable over $A$, then any Morley sequence in $p$ over $A$ is an indiscernible set over $A$. In particular, a Morley sequence in $p$ (over $\C$) is an indiscernible set over $\C$. 
\end{fact}

The following observation was made in \cite[Proposition 1.2]{11}.
\begin{fact}\label{dcl}
Let $p=tp(a/\C)$ be a type generically stable over $A$, and assume that 
$b \in dcl(A,a)$. Then $tp(b/\C)$ is also generically stable over $A$.
\end{fact}


The next lemma is a variant of a similar result for NIP groups (see \cite[Theorem 1.0.5]{Wa}), and its proof is very similar to the proof of \cite[Theorem 1.0.5]{Wa}. This is also a slight modification of \cite[Lemma 2.1(i)]{11}; in fact, it easily implies \cite[Lemma 2.1(i)]{11}.

\begin{lemma}\label{chain}
Let $G$ be a group which is $\emptyset$-definable in $\mathfrak{C}$ by a formula $G(x)$. 
Assume that $p\in S(\mathfrak{C})$
is generically stable over $A$. 
Let $H(x,\bar{z};y)$ be a formula defining a family of groups $H(G,c;g)$, $g\in (p|A)(\C)$, $c\in D$ ($D$ is a definable set).
Then there is $N<\omega$ such that for any $c\in D$, $n \in \omega$ and 
$(g_{1},\dots,g_{n})\models p^{(n)}|A$ there are $i_1,\dots,i_N \in \{ 1, \dots,n\}$ for which
\[\bigcap_{i=1}^{n}H(G,c;g_{i})=\bigcap_{j=1}^N H(G,c;g_{i_j}).\]
\end{lemma}
{\em Proof.}
Let $m>0$ be such as in the definition of generic stability for $p$ and $H(x,\bar{z} ;y)$. We will show that $N:=2m$ satisfies our requirements. Suppose it is not the case. Then there is $n>N$ (even $n=N+1$ works) such that for some $c \in D$ and $(g_{1},\dots,g_{n})\models p^{(n)}|A$ the intersection $\bigcap_{i=1}^{n}H(G,c;g_{i})$ is not an intersection of at most $n-1$ groups among $H(G,c;g_{i})$, $i=1,\dots,n$. Hence, for each $j \in \{1,\dots,n\}$ there exists 
$$a_j \in \bigcap_{i \ne j} H(G,c;g_{i}) \setminus \bigcap_{i=1}^{n}H(G,c;g_{i}).$$
Put $b=\prod_{j=1}^{m} a_j$. We see that 
$$b \in H(G,c;g_{i}) \iff i \in \{m+1,\dots,n\},$$
which contradicts the choice of $m$. \hfill $\blacksquare$\\

Recall that a subset of a group $G$ is said to be left generic if finitely many left translates of this set cover $G$. A formula $\varphi(x)$ is left generic if the set $\varphi(G)$ is left generic. A type is said to be left generic if every formula in it is left generic.

\begin{definition}
Let $G$ be a group definable in $\C$ by a formula $G(x)$.
$G$ has fsg (finitely satisfiable generics) if there is a global type $p$ containing $G(x)$ and a model $M\prec \C$, of cardinality less than the degree of saturation of $\mathfrak C$, such that for all $g$, $gp$ is finitely satisfiable in $M$.
\end{definition}

Let $G$ be a group $\emptyset$-definable in $\C$.
By $G^{00}$ we will denote the smallest type-definable subgroup of bounded index (if it exists). 
We do not know whether $G^{00}$ always exists when $T$ is $\omega$-categorical. Notice, however, that if it exists, then, being $\emptyset$-invariant, it must be $\emptyset$-definable and of finite index in $G$.
 
The following fact was proved in \cite[Section 4]{4}.

\begin{fact}\label{fsg}
Suppose $G$ has fsg,
witnessed by $p$. Then:\\
(i) a formula is left generic iff it is right generic (so we say that it is generic),\\
(ii) $p$ is generic,\\
(iii) the family of nongeneric sets forms an ideal, so any partial generic type can be extended to a global generic type,\\
(iv) $G^{00}$ exists, it is type-definable over $\emptyset$, and it is the stabilizer of any global generic type of $G$.
\end{fact}

Recall \cite[Proposition 0.26]{3}.

\begin{fact}\label{uniqueness}
Suppose $G$ has fsg and $G^{00}$ is definable. Then $G^{00}$ has a unique global generic type.
\end{fact}

The next definition was introduced in \cite[Section 6]{7}.

\begin{definition}
(i) Let $G$ be a group definable in $\C$. We say that $G$ is generically stable if it has fsg and some global generic type is generically stable.\\
(ii) Let $R$ be a ring definable in $\C$. We say that $R$ is generically stable if its additive group is generically stable.
\end{definition}

We say that a group [or ring] definable in a non-saturated model is generically stable if the group [or ring] defined by the same formula in a monster model is such. 

When we are talking about an $\omega$-categorical, generically stable group [or ring], we mean a generically stable group $G$ [or ring $R$] definable in a monster model $\C$ of an $\omega$-categorical theory. Replacing $\C$ by $G$ [or by $R$] equipped with the structure induced from $\C$, neither $\omega$-categoricity nor generic stability is lost. So, whenever we want to prove some algebraic properties of $G$ [or $R$], we can assume that $\C=G$ [or $\C=R$] (possibly with some extra structure).

We say that $G$ is connected if it does not have a proper, definable subgroup of finite index, and we will say that $G$ is absolutely connected if it does not have a proper, type-definable subgroup of bounded index (i.e. $G=G^{00}$).

Recall some basic notions from ring theory. In this paper, rings are associative, but they are not assumed to contain 1 or to be commutative. An element $r$ of a ring $R$ is nilpotent of nilexponent $n$ if $r^n=0$ and $n$ is the smallest number with this
property. The ring is nil [of nilexponent $n$] if every element is nilpotent [of nilexponent $\leq n$ and there is an element of nilexponent $n$]. The ring is nilpotent of class $n$ if $r_1\cdots r_n=0$ for all $r_1,\ldots,r_n\in R$ and $n$ is the smallest number with this property. An element $r$ is null if $rR=Rr=\{0\}$. The ring is null if all its elements are. 

The Jacobson radical of a ring $R$, denoted by $J(R)$, is the collection of all
elements of $R$ satisfying the formula $\phi(x)=\forall y\exists z(yx+z+zyx=0)$ (that is, it is the set of all elements which generate quasi-regular left ideals.). Equivalently, $J(R)$ is the intersection of all the maximal regular left [or right] ideals, where a left ideal $I$ is said to be regular if there is $a \in R$ such that $x-xa \in I$ for all $x \in R$ (notice that for rings with $1$ all ideals are regular). For any ring $R$, $J(R)$ is a two-sided ideal. We say that $R$ is semisimple if $J(R)=\{0\}$. $R/J(R)$ is always a semisimple ring. For details on Jacobson radical see \cite[Chapter 1]{14}. 

Recall that a ring $R$ is a subdirect product of rings $R_i,i\in I$, if there is a monomorphism of $R$ into $\prod_{i\in I}R_i$ whose image projects onto each $R_i$. 
The following fact is \cite[Corollary 1]{13}.

\begin{fact}\label{subd}
If $R$ is a semisimple, $\omega$-categorical ring, then $R$ is a subdirect product of complete matrix 
rings over finite field. Moreover, only finitely many different matrix rings occur as subdirect factors.
\end{fact}

By \cite[Lemma 1.3]{13} and \cite{15} we have:

\begin{fact}\label{nilp}
If $R$ is an $\omega$-categorical ring, then $J(R)$ is nilpotent.
\end{fact}

\section{$\pmb{\omega}$-categorical, generically stable rings}

This section is devoted to the proof Theorem \ref{Thm2} from the introduction. 
After a reduction to the situation when there is a unique global generic type, our proof splits into two cases depending on whether the generic type has non-nilpotent or nilpotent realizations. If they are non-nilpotent, the proof is a slight elaboration of the proof of \cite[Theorem 2.1]{5}, which is based on Facts \ref{subd} and \ref{nilp}. For the reader's convenience, we include most of the details. The argument in the nilpotent case is completely different; in particular, it uses a variant of Fact \ref{solvable} and some ideas from the proof of \cite[Theorem 2.1(i)]{12}. It will be noticed in the course of the proof that if the ring in question is commutative, then it is enough to consider the non-nilpotent case.

As to the nilpotent case, if we knew that the answer to \cite[Question 3.3]{11} is positive, i.e. if we knew that the generic stability of a type $p$ implies the generic stability of all its powers $p^{(n)}$, $n\geq 1$, then we could use Fact \ref{solvable} in the proof. Since \cite[Question 3.3]{11} remains open, we have to use a certain variant (in fact, strengthening) of Fact \ref{solvable}. Literally, it will be a strengthening of \cite[Theorem 2.3]{11}, obtained by the same proof as in \cite{11} modulo obvious modifications and
applications of an appropriate strengthening of \cite[Lemma 2.2]{11} which is described below.

\begin{lemma}\label{main lemma from the old paper}
Let $G$ be a $\emptyset$-definable group in a monster model $\C$ of an $\omega$-categorical theory. Assume that $G_{1}\unlhd G$ is infinite, $\emptyset$-definable, and characteristically simple in $(G,\C)$, i.e. it has no non-trivial, proper subgroup which is invariant under conjugations by the elements of $G$ and invariant under $Aut(\C)$.
Let $p\in S(\C)$ be a type generically stable over $\emptyset$. Suppose that for some $\emptyset$-definable function $f$ and Morley sequence $(g_{i})_{i<\omega}$ in $p$ over $\emptyset$, $f(g_0,\dots,g_{k-1}) \in G_1 \setminus \{ 0 \}$. Assume additionally that whenever $(h_i)_{i<k}$ is a Morley sequence in $p$ over some $g \in G$, then the conjugate $f(h_0,\dots,h_{k-1})^g$ equals $f(h_0',\dots,h_{k-1}')$ for some Morley sequence $(h_i')_{i<k}$ in $p$ over $g$.

Then $G_{1}$ is abelian.
\end{lemma}
{\em Sketch of proof.} 
Define $$H=\bigcap_{i_{1}<\dots <i_{k}}G_{G_1}(f(g_{i_{1}},\dots ,g_{i_{k}})).$$ 

It was shown in the course of the proof of \cite[Lemma 2.2]{11} that $H$ is invariant under $Aut(\C)$ (this follows from the generic stability of $p$ over $\emptyset$). 

Now, we will show that $H$ is normal in $G$. Consider any $g \in G$. Choose a Morley sequence $(h_i)_{i<\omega}$ in $p$ over $g$. By the invariance of $H$ under $Aut(\C)$ and the uniqueness of a Morley sequence in $p$ over $\emptyset$ up to the type, we have that
$$H=\bigcap_{i_{1}<\dots <i_{k}}C_{G_1}(f(h_{i_{1}},\dots,h_{i_{k}})).$$
Thus,
$$H^g=\bigcap_{i_{1}<\dots <i_{k}}C_{G_1}(f(h_{i_{1}},\dots,h_{i_{k}})^g).$$

By assumption, for every $i_1<\dots<i_k$ there are $h_{i_1}',\dots,h_{i_k}'$ 
forming a Morley sequence in $p$ over $g$ for which $f(h_{i_0},\dots,h_{i_{k-1}})^g=f(h_{i_0}',\dots,h_{i_{k-1}}')$. Once again by the invariance of $H$ under $Aut(\C)$, we get that $H \leq C_{G_1}(f(h_{i_0}',\dots,h_{i_{k-1}}'))$. Therefore, $H\leq H^g$. This shows that $H$ is normal in $G$.

Knowing that $H$ is invariant under $Aut(\C)$ and normal in $G$, the rest of the proof of \cite[Lemma 2.2]{11} works in our context, yielding the desired conclusion. \hfill $\blacksquare$\\

Having Lemma \ref{main lemma from the old paper}, in order to get the next theorem (strengthening \cite[Theorem 2.3]{11})), one should apply the proof of \cite[Theorem 2.3]{11} with several obvious modifications.

\begin{theorem}\label{main theorem from the old paper}
Let $G$ be a group $\emptyset$-definable in a monster model $\C$ of an $\omega$-categorical theory. Let $p \in S(\C)$ be a type generically stable over $\emptyset$. Suppose that for some $\emptyset$-definable function $f$ and Morley sequence $(g_i)_{i<\omega}$ in $p$ over $\emptyset$, $tp(f(g_0,\dots,g_{k-1})/\emptyset)$ is a generic type of $G$. Assume additionally that whenever $(h_i)_{i<k}$ is a Morley sequence in $p$ over $A,g$ (for some $A \subseteq \C$ and $g \in G$), then the conjugate $f(h_0,\dots,h_{k-1})^g$ equals $f(h_0',\dots,h_{k-1}')$ for some Morley sequence $(h_i')_{i<k}$ in $p$ over $A,g$.
Then $G$ is solvable-by-finite.
\end{theorem}

We would like to remark that while repeating the proof of \cite[Theorem 2.3]{11} to get Theorem \ref{main theorem from the old paper}, the function $f$ is replaced by some other functions, but always defined by means of $f$ in a purely group-theoretic way (e.g. by taking quotients or iterated commutators) and that is why the property of $f$ described in the last but one sentence of Theorem \ref{main theorem from the old paper} is never lost (and we can use Lemma \ref{main lemma from the old paper}).


\begin{remark}\label{R^00}
Suppose $R$ is a ring definable in a monster model and that $R^{00}$ in the additive sense exists. Then $R^{00}$ is an ideal of $R$.
\end{remark}
{\em Proof.} Consider any $r \in R$. Let $f : R \to R$ be an additive homomorphism defined by $f(x)=rx$. Then $f[R^{00}]$ is a type-definable subgroup of $R$.

We need to show that $f[R^{00}] \subseteq R^{00}$. Suppose this is not the case. Then $A:=R^{00}\cap f[R^{00}]$ is 
a proper, type-definable subgroup of $f[R^{00}]$ of bounded index. Hence, $R^{00} \cap f^{-1}[A]$ is a proper, type-definable, bounded index subgroup of $R^{00}$, which is not possible. \hfill $\blacksquare$\\


\noindent
{\em Proof of Theorem \ref{Thm2}.}
Let $R$ be a generically stable ring which is definable in a monster model of an $\omega$-categorical theory. 
By Fact \ref{fsg}(iv), $R^{00}$ (in the additive sense) exists, and since it is $\emptyset$-invariant, it follows by $\omega$-categoricity that it is $\emptyset$-definable. So, it has finite index in $R$, and, by Remark \ref{R^00}, we can assume that $R=R^{00}$ 
(because, being an additive translate of a generically stable generic type, the generic type of $R^{00}$ is also generically stable by Fact \ref{dcl}). Fact \ref{nilp} tells us that 
$J(R)$ is nilpotent, so we can assume that $R$ is semisimple
(replacing $R$ by $R/J(R)$ and using Fact \ref{dcl}). 
If $R$ is finite, we are done, so we can assume that
$R$ is infinite.  By Fact \ref{uniqueness}, $R$ has a unique (global) generic type $p \in S_1(R)$. 
Thus, $p$ is generically stable over $\emptyset$. Without loss of generality, the monster model is just the ring $R$ (possibly with some extra structure).

As it was mentioned at the beginning of this section, the proof splits into two cases depending on whether the realizations of $p$ are non-nilpotent or nilpotent. Notice, however, that in the special case of commutative rings, all non-zero elements of $R$ are non-nilpotent 
(since, being semisimple, $R$ does not have non-trivial nil ideals). So, for commutative rings it is enough to consider only the first case.\\[2mm]
{\bf Case 1.} Realizations of $p$ are not nilpotent. \\[2mm]
By Fact \ref{subd}, we can assume that $R$ is a subring of $\prod_{i\in I}R_{i}$, where each $R_{i}$ is finite, and 
$|\{R_{i}:i\in I\}|<\omega$. Let $\pi_{i}$ be the projection onto the $i$-th coordinate.
For $i_{0},\dots,i_{n}\in I$ and $r_{0}\in R_{i_{0}},\dots,r_{n}\in R_{i_{n}}$, we define
$$R_{i_{0},\dots,i_{n}}^{r_{0},\dots,r_{n}}=\left\{r\in R : \bigwedge_{j=0}^{n}\pi_{i_{j}}(r)=r_{j}\right\}.$$

\begin{claim}
There are $i_{0},i_{1},\dots\in I$, non-nilpotent elements $r_{j}\in R_{i_{j}}$ and a Morley 
sequence $(\eta_{i})_{i<\omega}$ in $p$ over $\emptyset$ such that $\eta_{n}\in R_{i_{0},\dots,i_{n-1},i_{n}}^{0,\dots,0,r_{n}}$ for every $n<\omega$. 
\end{claim}

\noindent
{\em Proof of Claim 1.} Assume that we have already constructed $i_{0},\dots,i_{n}$, $r_{0},\dots,r_{n}$ and $(\eta_{i})_{i\leq n}$. 
Let $p_n=p|{(\eta_{i})_{i\leq n}} \in S_1((\eta_i)_{i \leq n})$.
If $R_{i_{0},\dots,i_{n}}^{0,\dots,0}\cap p_n (R)=\emptyset$, then 
$R\setminus p_n(R)$ is generic in $R$ (since $R_{i_{0},\dots,i_{n}}^{0,\dots,0}$ has finite index in $R$), 
so, by compactness, there is $\phi\in p_n$ such that $\neg\phi$ is generic. 
From Fact \ref{fsg}(iii), we get that $\{\neg\phi\}$ extends to a global generic type, which contradicts the uniqueness of the generic type in $R$. So, $R_{i_{0},\dots,i_{n}}^{0,\dots,0}\cap p_n(R)\neq\emptyset$.
Take $\eta_{n+1}\in R_{i_{0},\dots,i_{n}}^{0,\dots,0}\cap p_n(R)$.
By the assumption of Case 1, $\eta_{n+1}$ is non-nilpotent. Since $|\{R_{i}:i\in I\}|<\omega$, there is $i_{n+1}\in I$ such that $\pi_{i_{n+1}}(\eta_{n+1})$ is non-nilpotent. We put $r_{n+1}:=\pi_{i_{n+1}}(\eta_{n+1})$.
So, we have found $i_{n+1}$, $r_{n+1}$ and $\eta_{n+1}$ with the desired properties.\hfill $\square$\\ 


By $\omega$-categoricity, the two-sided ideals $RrR$, $r \in R$, are uniformly definable (because $\omega$-categoricity implies that there is $K$ such that every element of any $RrR$ is the sum of at most $K$ elements of the form $r_1rr_2$ for $r_1,r_2 \in R\cup \{ 1\}$). Thus, there is a formula $H(x,z;y)$ expressing that $x\in R(z-y)R$. Let $N$ be as in Lemma \ref{chain} for the type $p$, formula $H(x,z;y)$ and $D:=R$.


\begin{claim}
There are natural numbers $n(0)<n(1)<\dots<n(N)$ 
and a Morley sequence $(a_{i})_{i\leq N}$ (of length $N+1$) in $p$ over $\emptyset$ 
such that 
\begin{equation*}\tag{$*$}
a_0\in R_{i_{n(0)},\dots,i_{n(N)}}^{r_{n(0)},0,\dots,0},a_{1}\in R_{i_{n(0)},\dots,i_{n(N)}}^{0,r_{n(1)},0,\dots,0}, \dots ,a_N\in R_{i_{n(0)},\dots,i_{n(N)}}^{0,\dots,0,r_{n(N)}}.
\end{equation*} 
\end{claim}
{\em Proof of Claim 2.} 
First, we will find natural numbers $$n(0)<n'(0)<n(1)<n'(1)<\dots <n(N-1)<n'(N-1)<n(N)$$ such that for $a_{k}:=\eta_{n(k)}-\eta_{n'(k)}$, $k=0,\dots,N-1$, and $a_N:=\eta_{n(N)}$ the condition $(*)$ is satisfied. This follows exactly as in the proof of Claim 2 in \cite[Theorem 2.1]{5}, but we give some details for completeness. 

Let $c=\max_{i\in I}|R_{i}|$. Define recursively numbers $L_N,\dots,L_1,L_0$:
$$
\begin{array}{lll}
L_N & = & c+1,\\
L_{N-k} & = & c^{L_N+\dots+L_{N-k+1}+1}+1 \;\;\mbox{for}\;\; k=1,\dots,N-1,\\
L_0& = & 0.
\end{array}
$$
%
Put $I_N = \{ L_0+\dots +L_N\}$, and define intervals $I_0,\dots,I_{N-1}$ as
$$I_k=[L_0+\dots+L_k,L_0+\dots+L_{k+1}-1].$$
For each $k\in \{0,\dots,N-1\}$, by the pigeonhole principle, we can find two natural numbers $n(k)<n'(k)$ in $I_{k}$ such that $\pi_{i_{j}}(\eta_{n'(k)})=\pi_{i_{j}}(\eta_{n(k)})$ for every $j\in I_{k+1}\cup\dots\cup I_{N}$. Put additionally $n(N)=L_0+\dots+L_N$.
Now, it is easy to check that $(*)$ is satisfied for $a_{k}:=\eta_{n(k)}-\eta_{n'(k)}$, $k=0,\dots,N-1$, and $a_N:= \eta_{n(N)}$. 

It remains to show that $(a_{i})_{i\leq N}$ is a Morley sequence in $p$ over $\emptyset$.
Fix any $k<N$. We have that $\eta_{n'(k)}\models p|(\eta_{n'(i)})_{i<k}(\eta_{n(i)})_{i\leq k}$. By the uniqueness of the generic type in $R$, we get that $\eta_{n(k)}-\eta_{n'(k)}\models p|(\eta_{n'(i)})_{i<k}(\eta_{n(i)})_{i< k}$, so $a_{k}\models p|(a_{i})_{i<k}$. It is also clear that $a_N \models p|(a_i)_{i<N}$. This shows that $(a_{i})_{i\leq N}$ is a Morley sequence in $p$ over $\emptyset$. \hfill $\square$\\

Let $c=\sum_{i\leq N}a_{i}$ and $b_{j}=\sum_{i\neq j}a_{i}=c-a_{j}$ for $j=0,\dots,N$.
Using Claim 2 and the choice of $N$, we reach a final contradiction in the same way as in the proof of \cite[Theorem 2.1]{5}. Namely, 

\begin{equation*}\tag{$**$}
\pi_{i_{n(j)}}[Rb_0R\cap\dots \cap Rb_NR]= \{0\}\;\, \mbox{for}\;\, j=0,\dots, N. 
\end{equation*}
On the other hand, $\prod_{k \ne j}b_k \in \bigcap_{k \ne j}Rb_kR$ for $j=0,\dots,N$. We also have that $\pi_{i_{n(j)}}[\prod_{k \ne j}b_k]=r_{n(j)}^N \ne 0$ as $r_{n(j)}$ is non-nilpotent. So, 
\begin{equation*}\tag{$***$}
 \pi_{i_{n(j)}}\left[\bigcap_{k \ne j}Rb_kR\right] \ne \{ 0 \} \;\, \mbox{for}\;\, j=0,\dots, N.
\end{equation*}
 By $(**)$ and $(***)$, $Rb_0R\cap \dots\cap Rb_N R \ne \bigcap_{k\ne j}Rb_kR$ for all $j=0,\dots,N$. This is a contradiction with the choice of $N$, because $Rb_iR=R(c-a_i)R=H(R,c;a_i)$ and $(a_i)_{i\leq N}$ is a Morley sequence in $p$ over $\emptyset$.\\[2mm]
%
%
{\bf Case 2.} Realizations of $p$ are nilpotent. \\[2mm]
By $\omega$-categoricity, $R$ has a finite characteristic $c$. Put $R_1=R\times \Z_c$, and define $+$ and $\cdot$ on $R_1$ by $(a,k)+(b,l)=(a+b,k+_cl)$ and $(a,k)\cdot(b,l)=(ab+l\times a+k\times b,k\cdot_cl)$, where $+_c$ and $\cdot_c$ are addition and multiplication modulo $c$, and $l\times a:=a+\dots+a$ ($l$-many times).
Then $R_1$ is a ring with $1$ interpretable in $R$, and $R$ is a two-sided ideal of finite index in $R_1$.
Let $G$ be a subgroup of $\GL_3(R_1)$ generated by $\{t_{ij}(\alpha):\alpha \in R,\; i,j\in \{1,2,3\},\; i\neq j \}\cup\{t_{j}(\beta):\beta \in (1+R)\cap R_1^*, \; j\in \{1,2,3\} \}$, where $t_{ij}(\alpha)$ is the matrix with $1$'s on the diagonal, $\alpha$ on the $(i,j)$-th position and $0$'s elsewhere, and 
$t_j(\beta)$ is the matrix with $\beta$ on the $(j,j)$-th position, $1$'s on the rest of the diagonal and $0$'s elsewhere. 
Since $G$ is invariant over finitely many parameters (over which $R_1$ is defined), it follows by $\omega$-categoricity that it is definable.
Let $(a_{ij})_{1\leq i,j\leq 3}\models p^{(9)}$ (in a bigger monster model $\C\succ R$); note that the order of $a_{ij}$'s is irrelevant, because $p$ is generically stable and we have Fact \ref{indiscernible}. Define 
$$A=\left( \begin{array}{ccc}
1+a_{11} & a_{12} & a_{13} \\
a_{21} & 1+a_{22} & a_{23} \\
a_{31} & a_{32} & 1+a_{33} \end{array} \right). $$

\begin{claim}
$A\in G(\C)$, where $G(\C)$ is the interpretation of $G$ in $\C$.
\end{claim}
{\em Proof of Claim 3.}
The idea is to show that we can transform $A$ to the identity matrix by Gaussian elimination process in which all elementary matrices belong to $G(\C)$ (because then $BA=I$ for some $B \in G(\C)$, so $A=B^{-1} \in G(\C)$). 

The following well-known remark is fundamental for our process: if $r \in R$ satisfies $r^n=0$ for some $n$, then $(1+r)(1-r+r^2-\dots \pm r^{n-1})=1$, so $(1+r)^{-1}\in (1+R\cap dcl(r))\cap R_1^*$. 

Now, we describe the first step of the process. We have 

$$\begin{array}{cc}
t_{21}\left(-a_{21}(1+a_{11})^{-1}\right)A=\\ 
\left( \begin{array}{ccc}
1+a_{11} & a_{12} & a_{13} \\
0 & 1+a_{22}-a_{21}(1+a_{11})^{-1}a_{12}& a_{23}-a_{21}(1+a_{11})^{-1}a_{13}\\
a_{31} & a_{32} & 1+a_{33} \end{array} \right). \end{array} $$ 
But $a_{21}(1+a_{11})^{-1}a_{12}\in \C\cap dcl((a_{ij})_{(i,j)\neq (2,2)})$, so, by the uniqueness of the generic type, $$b_{22}:=a_{22}-a_{21}(1+a_{11})^{-1}a_{12}\models p|R,(a_{ij})_{(i,j)\neq (2,2)}.$$ Similarly, $$b_{23}:=a_{23}-a_{21}(1+a_{11})^{-1}a_{13}\models p|R,(a_{ij})_{(i,j)\neq (2,3)}.$$ 
Therefore, $$t_{21}\left(-a_{21}(1+a_{11})^{-1}\right)A=\left( \begin{array}{ccc}
1+a_{11} & a_{12} & a_{13} \\
0 & 1+b_{22} & b_{23} \\
a_{31} & a_{32} & 1+a_{33} \end{array} \right),$$
where $((a_{ij})_{(i,j) \ne (2,2),(2,3)}, b_{22},b_{23}) \models p^{(8)}$.
Continuing Gaussian elimination in this way, we obtain a matrix $C \in G(\C)$ such that
$$CA= \left( \begin{array}{ccc}
1+b_1 & 0& 0 \\
0 & 1+b_2 & 0 \\
0 & 0 & 1+b_3 \end{array} \right) =t_1(1+b_1)t_2(1+b_2)t_3(1+b_3)$$ for some $(b_1,b_2,b_3)\models p^{(3)}$. But $t_1(1+b_1),t_2(1+b_2),t_3(1+b_3)\in G(\C)$, so we conclude that $A\in G(\C)$. \hfill $\square$
%
\begin{claim}
All translates of the type $q:=tp(A/R)$ by the elements of $G$ are finitely satisfiable in some small model. Thus, $G$ has fsg and $q$ is a global generic type of $G$.
\end{claim}
{\em Proof of Claim 4.}
We will show that every translate of $A$ by an element of $G$ belongs to the set
$$Z:=\{ (k_{ij}+b_{ij})_{i,j\in \{1,2,3\}} :k_{ij}\in \Z_c,\; (b_{ij})_{i,j\in \{1,2,3\}}\models p^{(9)} \}.$$ 
This will complete the proof of Claim 4, as every element of $Z$ is in the definable closure of $\Z_c$ and some realization of the type $p^{(9)}$, and $p^{(9)}$ is finitely satisfiable in some small model. So, it suffices to show that $Z$ is invariant under multiplication by the elements of the set 
$\{t_{ij}(\alpha):\alpha \in R,\; i,j\in \{1,2,3\},\; i\neq j \}\cup\{t_{j}(\beta):\beta \in (1+R)\cap R_1^*,\; j\in \{1,2,3\} \}$ (notice that this set is closed under the group inversion). Choose any $B=(k_{ij}+b_{ij})_{i,j\in \{1,2,3\}} \in Z$. 

First, consider any element of $Z$ of the form $t_{j}(\beta)$, where $\beta=1+r\in R_1^*$ for some $r \in R$. 
Denote the coefficients of the matrix $t_{j}(\beta)B$ by $d_{im}$ ($i,m\in \{1,2,3\}$). 
Then $d_{im}=k_{im}+b_{im}$ for all $m$ and $i\neq j$. Take any $m\in \{1,2,3\}$. Then $d_{jm}=\beta (k_{jm}+b_{jm})=(1+r)(k_{jm}+b_{jm})=k_{jm}+k_{jm}\times r+(1+r)b_{jm}.$ 
Since multiplication by $(1+r)$ is a definable automorphism of $(R,+)$, by the uniqueness of the generic type, we get that $(1+r)b_{jm}\models p|R,(b_{il})_{(i,l)\neq (j,m)}$, and hence $k_{jm}\times r+(1+r)b_{jm} \models p| R,(b_{il})_{(i,l)\neq (j,m)}$. This easily implies that $t_{j}(\beta)B\in Z$. 

Now, consider any $t_{ij}(\alpha)$, where $i \ne j$ and $\alpha \in R$. Denote the coefficients of $t_{ij}(\alpha)B$ by $f_{ij}$ ($i,j\in \{1,2,3\}$). Choose any $m\in \{1,2,3\}$. For all $l\neq i$ we have $f_{lm}=k_{lm}+b_{lm}$. Moreover, $f_{im}=k_{im}+b_{im}+\alpha (k_{jm}+b_{jm})$. But $\alpha (k_{jm}+b_{jm})\in dcl(R,b_{jm})$, so $b_{im}+\alpha (k_{jm}+b_{jm})\models p|R,(b_{pq})_{(p,q)\neq (i,m)}$. Hence, $t_{ij}(\alpha)B\in Z$. This completes the proof of Claim 4. \hfill $\square$\\

If we knew that $p^{(9)}$ is generically stable (recall that we know that $p$ is generically stable), then Claim 4 would imply that $G$ is generically stable, so, by Fact \ref{solvable}, we would conclude that $G$ is solvable-by-finite and we could turn to the last paragraph of the proof. Since we do not know whether $p^{(9)}$ is generically stable, we will prove Claim 5 below and then apply Theorem \ref{main theorem from the old paper} in order to get that $G$ is solvable-by-finite.
 
Adding to the language the appropriate parameters, we can assume that everything is definable over $\emptyset$.

\begin{claim}
Let a $\emptyset$-definable function $f: M_{3 \times 3}(R) \to M_{3\times 3}(R_1)$ be defined by 
$$f((x_{ij})_{1\leq i,j\leq 3}))=\left( \begin{array}{ccc}
1+x_{11} & x_{12} & x_{13} \\
x_{21} & 1+x_{22} & x_{23} \\
x_{31} & x_{32} & 1+x_{33} \end{array} \right).$$
Then, whenever $(h_{ij})_{1\leq i,j\leq3}$ is a Morley sequence in $p$ over $A,g$ (for some $A\subseteq R$ and $g \in G$), then
$f((h_{ij})_{1\leq i,j \leq 3}) \in G$ and $f((h_{ij})_{1\leq i,j \leq 3})^g=f((h_{ij}')_{1\leq i,j\leq3})$ for some Morley sequence $(h_{ij}')_{1\leq i,j\leq3}$ in $p$ over $A,g$.
\end{claim}
{\em Proof of Claim 5.} Let $(h_{ij})_{1\leq i,j\leq3}$ be a Morley sequence in $p$ over $A,g$. The fact that $f((h_{ij})_{1\leq i,j \leq 3}) \in G$ follows from Claim 3 and the uniqueness of a Morley sequence in $p$ over $\emptyset$ up to the type.

For the second part, first notice that it is enough to prove the statement for $g$ of the form
$t_{ij}(\alpha)$ (for $\alpha \in R$ and distinct $i,j\in \{1,2,3\}$) and $t_{j}(\beta)$ (for $\beta \in (1+R)\cap R_1^*$ and $j\in \{1,2,3\})$, which follows easily from the fact that $G$ is generated by these elements and by the independence of the choice of a Morley sequence $(h_{ij})_{1\leq i,j\leq3}$ in $p$ over $A,g$. Next, apply a similar argument to the proof of Claim 4 in order to get that the conjugates by the elements of the form $t_{ij}(\alpha)$ or $t_j(\beta)$ have the desired property. \hfill $\square$\\

From Claims 4 and 5 and Theorem \ref{main theorem from the old paper}, we conclude that $G$ is solvable-by-finite.

The rest of the proof follows exactly as in \cite[Theorem 2.1(i)]{12}, but we will repeat the argument for the reader's convenience.
Let $H$ be a normal subgroup of $G$ of finite index, which is solvable. We have the following well-known formulas:\\ 
\begin{equation}\tag{$\dag$}
t_{ij}(\alpha)t_{ij}(\beta)=t_{ij}(\alpha +\beta) \;\; \mbox{and}\;\; [t_{ik}(\alpha),t_{kj}(\beta)]=t_{ij}(\alpha\beta)
\end{equation}
for pairwaise distinct $i,j,k$. Define 
$I=\{r\in R:(\forall i\neq j) t_{ij}(r)\in H\}$. Using the normality of $H$ in $G$ and $(\dag)$, we see that $I\vartriangleleft R$. 
If $|R/I|\geq \omega$, then, by Ramsey theorem, 
for some distinct $i,j\in \{1,2,3\}$ there are $r_k$, $k<\omega$, such that $t_{ij}(r_n-r_m)\notin H$ for every $n<m<\omega$. 
But, by $(\dag)$, $t_{ij}(r_n-r_m)=t_{ij}(r_n)t_{ij}(r_m)^{-1}$, which contradicts the finiteness of $[G:H]$.
So, $|R/I|<\omega$. Since $H$ is solvable, there exists $n$ for which the $n$-th derived subgroup $H^{(n)}$ is trivial.
 Then $(\dag)$ implies that for every $r_{1},\dots,r_{2^n}\in I$ and distinct $i,j\in \{1,2,3\}$ we have $t_{ij}(r_{1}\dots r_{2^n})\in H^{(n)}=\{e\}$, so $r_{1}\dots r_{2^n}=0$.
This shows that $I$ is a nilpotent ideal of $R$ of finite index. \hfill $\blacksquare$

\section{$\pmb{\omega}$-categorical, generically stable groups}

The goal of this section is to prove Theorem \ref{Thm1} from the introduction. In the final part of the proof, we will apply Fact \ref{solvable} and the argument from page 490 of \cite{6}. However, in order to do that, first we need to prove a certain non-trivial lemma (a variant of \cite[Corollary 3.5]{6}), which uses some ideas from the proof of Theorem \ref{Thm2} in Case 1 and from the final part of the proof of \cite[Corollary 3.17]{12}. 

It is worth recalling that \cite[Theorem 3.15]{12} says that whenever each ring interpretable in a given $\omega$-categorical structure (in which all definable groups have connected components) is nilpotent-by-finite, then each solvable group definable in this structure is also nilpotent-by-finite. This was used in \cite{5} to conclude that $\omega$-categorical groups with NIP and fsg are nilpotent-by-finite (using the fact that they are solvable-by-finite and that each $\omega$-categorical ring with NIP is nilpotent-by finite). The reason why, having Theorem \ref{Thm2}, we cannot apply Fact \ref{solvable} and {\cite[Theorem 3.15]{12} in order to get Theorem \ref{Thm1} is that rings interpretable in a given $\omega$-categorical, generically stable group need not to be generically stable. In the proof of Lemma \ref{main lemma} below, we undertake a detailed analysis of the relevant interpretable rings.

Recall a few basic definitions. Let $H$ be a group. The commutator of $h_0,h_1 \in H$ is defined as $[h_0,h_1]=h_0^{-1}h_1^{-1}h_0h_1$.
The iterated commutators $\gamma_n$ on a group $H$ are defined inductively as follows: 
$$\gamma_1(h_0)=h_0\;\; \mbox{and}\;\; \gamma_{n+1}(h_0,\dots,h_n)=[\gamma_n(h_0,\dots,h_{n-1}),h_n].$$
The lower central series $H=\Gamma_1(H) \geq \Gamma_2(H) \geq \dots$ of $H$ is defined by:
$$\Gamma_1(H)=H\;\; \mbox{and}\;\; \Gamma_{n+1}(H) = [\Gamma_n(H),H].$$ 
%
 The following formulas for commutators will be very useful:
\begin{equation}
[x,zy]=[x,y][x,z]^y\;\; \mbox{and}\;\; [xz,y]=[x,y]^z[z,y]. 
\end{equation}
Recall that, in this paper, a group $H$ definable in a monster model is said to be absolutely connected if $H^{00}=H$.

\begin{lemma}\label{commutators}

Let $H$ be an absolutely connected group with fsg ($H$ is definable over $A$ in a monster model). Then for every $n\in \omega \setminus \{ 0 \}$ each element of $\Gamma_n(H)$ is a product of conjugates of elements from the set $\{\gamma_n(g_0,\dots,g_{n-1}):(g_0,\dots,g_{n-1})\models p^{(n)}|A\}$, where $p$ is the unique global generic type in $H$.
\end{lemma}
{\em Proof.} 
Using (1), we easily get by induction that for every $x_0,\dots,x_k,y$ the commutator $[x_0\dots x_k,y]$ is a product of conjugates of commutators $[x_i,y],i=0,\dots,k$.
To prove the lemma, we proceed by induction on $n$. The induction starts, since every element of $H$ is a product of two realizations of $p|A$ (which follows from the uniqueness of the generic type). Suppose that the conclusion of the lemma is satisfied for $n$. Take any $g\in \Gamma_{n+1}(H)$. Then $g=[a,b]$ for some $a\in \Gamma_n(H)$ and $b\in H$. By the inductive hypothesis, $a=\prod_{i<l}\gamma_n(\overline{g_i})^{c_i}$ for some $\overline{g_i}\models p^{(n)}|A$ and $c_i\in H$ (for $i = 0,\dots,l-1$). So, by (1), $[a,b]=\prod_{i<l}[\gamma_n(\overline{g_i})^{c_i},b]^{d_i}$ for some $d_i\in H$, $i=0,\dots,l-1$. Choose $b_1,b_2\models p|A,(\overline{g_{i}},c_i)_{i<l}$ such that $b_1b_2=b$. Then, by (1), for every $i<l$ we have

$$\begin{array}{llllll}
[\gamma_n(\overline{g_i})^{c_i},b]^{d_i} &= & [\gamma_n(\overline{g_i})^{c_i},b_2]^{d_i}[\gamma_n(\overline{g_i})^{c_i},b_1]^{b_2d_i} \\ 
& = & [\gamma_n(\overline{g_i}),b_2^{c_i^{-1}}]^{c_id_i}[\gamma_n(\overline{g_i}),b_1^{c_i^{-1}}]^{c_ib_2d_i}\\ 
& = &\gamma_{n+1}(\overline{g_i},b_2^{c_i^{-1}})^{c_id_i}\gamma_{n+1}(\overline{g_i},b_1^{c_i^{-1}})^{c_ib_2d_i}.
\end{array}$$ 
By the uniqueness of the generic type, $b_1^{c_i^{-1}},b_2^{c_i^{-1}}\models p| A,\overline{g_i}$. So, $(\overline{g_i}, b_1^{c_i^{-1}}),(\overline{g_i},b_2^{c_i^{-1}})\models p^{(n+1)}|A$. This completes the proof of the lemma. \hfill $\blacksquare$

\setcounter{claim}{0}

\begin{lemma}\label{main lemma}
We work in a monster model $\C$ of an $\omega$-categorical theory.
Let $H$ be a nilpotent, absolutely connected, generically stable group definable in $\C$ acting definably and by automorphisms on a definable vector space $V$ over $F:=GF(q^a)$ ($q$ is a prime number), and assume that $H$ has no elements of order $q$. Then $C_H(V)=H$. 
\end{lemma}
{\em Proof.} 
 First, we will prove the following claim.
\begin{claim} If $H$ is a nilpotent group of class $n$ which satisfies the assumptions of the lemma, then $[\Gamma_n(H):C_{\Gamma_n(H)}(V)]<\omega$. \end{claim}
{\em Proof of Claim 1.} 
We can assume that everything is $\emptyset$-definable in $\C$. Then the unique global generic type $p$ of $H$ is generically stable over $\emptyset$.
Put $A=\Gamma_n(H)$. 
As in \cite{6}, we define $W$ as the sum of all finite dimensional $FA$-submodules of $V$. 
Then $W$ is a definable in $\C$ and invariant under $A$ subspace of $V$. By \cite[Proposition 3.4]{6}, it is enough to show that $W=V$. Suppose for a contradiction that $W\subsetneq V$, and put $\overline{V}=V/W$.
Exactly as in in the proof of \cite[Corollary 3.5]{6}, we get:
\begin{equation}\tag{$*$}
\mbox{The $FA$-module $\overline{V}$ has no non-trivial, finite dimensional $FA$-submodules.}
\end{equation}

Choose a non-trivial $v\in \overline{V}$, and put $V_0=Lin_F(Av)$. By $\omega$-categoricity, $V_0$ is interpretable. Let $R$ be the ring of endomorphisms of $V_0$ generated by $A$. 
Since $A$ is commutative, $R$ is a commutative ring interpretable in $\C$. 
Adding some parameters to the language, we can assume that $R$ is interpretable in $\C$ over $\emptyset$. 


Let $(g_i)_{i<\omega}\models p^{(\omega)}|\emptyset$.
We will show that
\begin{equation}\tag{$**$}
\gamma_n(g_0,\dots,g_{n-2},g_n)-\gamma_n(g_0,\dots,g_{n-2},g_{n-1})\in J(R). 
\end{equation}
Suppose for a contradiction that this is not the case. 
Put $$a=\gamma_{n-1}(g_0,\dots,g_{n-2}).$$

By Fact \ref{subd}, we can assume that $R/J(R)$ is a subring of $\prod_{i\in I}R_{i}$, where each $R_{i}$ is finite and 
$|\{R_{i}:i\in I\}|<\omega$. Let $\pi_{i}:R\to R_i$ be the quotient map $R\to R/J(R)$ composed with the projection onto the $i$-th coordinate.
For $i_{0},\dots,i_{m}\in I$ and $r_{j}\in R_{i_{j}}$, we define 
$$R_{i_{0},\dots,i_{m}}^{r_{0},\dots,r_{m}}=\left\{r\in R : \bigwedge_{j=0}^{m}\pi_{i_{j}}(r)=r_{j}\right\}.$$ 

\begin{subclaim}
There are $i_{0},i_{1},\dots\in I$, non-nilpotent elements $r_{j}\in R_{i_{j}}$ and a Morley 
sequence $(\eta_{i})_{i<\omega}$ in $p$ over $a$ such that $[a,\eta_{2m+1}]-[a,\eta_{2m}]\in R_{i_{0},\dots,i_{m-1},i_{m}}^{0,\dots,0,r_{m}}$ for every $m<\omega$.
\end{subclaim}
{\em Proof of Subclaim 1.} 
Suppose we have already constructed $i_{0},\dots,i_{k-1}$, $r_{0},\dots,r_{k-1}$ and $(\eta_{i})_{i<2k}$. 
Let $(h_{i})_{i<\omega}\models p^{(\omega)} |a,(\eta_{i})_{i<2k}$. Choose $j<l<\omega$ such that $\pi_{i_m}([h_{j},a])=\pi_{i_m}([h_{l},a])$ for every $m< k$. Put $\eta_{2k}=h_{j}$ and $\eta_{2k+1}=h_{l}$. Since $R/J(R)$ is a semisimple, commutative ring, the only nilpotent element of $R/J(R)$ is zero. Hence, by our assumption that $(**)$ does not hold, $([a,\eta_{2k+1}]-[a,\eta_{2k}])/J(R)$ is non-nilpotent. Therefore, since $|\{R_{i}:i\in I\}|<\omega$, there is $i_{k}\in I$ such that $\pi_{i_{k}}([a,\eta_{2k+1}]-[a,\eta_{2k}])$ is non-nilpotent. 
Putting $r_{k}=\pi_{i_{k}}([a,\eta_{2k+1}]-[a,\eta_{2k}])$, the construction is completed.\hfill $\square$\\

Now, using the assumption that $H$ is nilpotent of class $n$ and (1), we get 
$$\begin{array}{lllll}
([a,\eta_{2k+1}]-[a,\eta_{2k}])[a,\eta_{2k}^{-1}] & = & [a,\eta_{2k+1}][a,\eta_{2k}^{-1}]-[a,\eta_{2k}][a,\eta_{2k}^{-1}]\\
& = & [a,\eta_{2k+1}][a,\eta_{2k}^{-1}]^{\eta_{2k+1}}-[a,\eta_{2k}][a,\eta_{2k}^{-1}]^{\eta_{2k}} \\ 
&= & [a,\eta_{2k}^{-1}\eta_{2k+1}]-1
\end{array}$$ for any $k<\omega$.
But $[a,\eta_{2k}^{-1}]$ is an invertible element of $R$, so putting $h_k=\eta_{2k}^{-1}\eta_{2k+1}$ for $k<\omega$, we obtain that $[a,h_k]-1\in R_{i_{0},\dots,i_{k-1},i_{k}}^{0,\dots,0,s_{k}}$ for every $k<\omega$, where $s_{k}\in R_{i_{k}}$ are non-nilpotent. Also, by the uniqueness of the generic type, $(h_i)_{i<\omega}$ is a Morley sequence in $p$ over $a$. 

Let $H(x,\overline{z};y)$ be a formula expressing that $x\in R(z_1-[z_2,y]+1)R$, where $\overline{z}=(z_1,z_2)$. Choose $N$ as in Lemma \ref{chain} for the type $p$, formula $H(x,\overline{z};y)$ and $D:=R \times \{ a\}$. 

Now, exactly as in Claim 2 in the proof of Theorem \ref{Thm2}, we find $$n(0)<n'(0)<n(1)<n'(1)<\dots <n(N-1)<n'(N-1)<n(N)<n'(N)$$ such that for $a_{k}:=([a,h_{n(k)}]-1)-([a,h_{n'(k)}]-1)$, $k=0,\dots,N$, we have: $$a_{0}\in R_{i_{n(0)},\dots,i_{n(N)}}^{s_{n(0)},0,\dots,0},a_{1}\in R_{i_{n(0)},\dots,i_{n(N)}}^{0,s_{n(1)},0,\dots,0}, \dots ,a_{N}\in R_{i_{n(0)},\dots,i_{n(N)}}^{0,\dots,0,s_{n(N)}}.$$
By the fact that $H$ is nilpotent of class $n$ and (1), for all $k=0,\dots,N$ we have 
$$a_k':=a_{k}[a,h_{n'(k)}^{-1}]=[a,h_{n(k)}][a,h_{n'(k)}^{-1}]-[a,h_{n'(k)}][a,h_{n'(k)}^{-1}]=[a,h_{n'(k)}^{-1}h_{n(k)}]-1.$$ 
Since each $[a,h_{n'(k)}^{-1}]$ is invertible in $R$, we obtain
$$a_{0}'\in R_{i_{n(0)},\dots,i_{n(N)}}^{s_{n(0)}',0,\dots,0},a_{1}'\in R_{i_{n(0)},\dots,i_{n(N)}}^{0,s_{n(1)}',0,\dots,0}, \dots ,a_{N}'\in R_{i_{n(0)},\dots,i_{n(N)}}^{0,\dots,0,s_{n(N)}'},$$
for some non-nilpotent elements $s_{n(0)}' \in R_{i_{n(0)}},\dots, s_{n(N)}' \in R_{i_{n(N)}}$.
Moreover, we see that $(h_{n'(k)}^{-1}h_{n(k)})_{k<\omega}$ is a Morley sequence in $p$ over $a$. So, by the choice of $N$,
the argument after Claim 2 in the proof of Theorem \ref{Thm2} leads to a contradiction, which
completes the proof of $(**)$.

By $(**)$ together with the fact that $H$ is nilpotent of class $n$ and (1), we obtain 
$$\begin{array}{lll}
J(R)&\ni&(\gamma_n(g_0,\dots,g_{n-2},g_n)-\gamma_n(g_0,\dots,g_{n-2},g_{n-1}))\gamma_n(g_0,\dots,g_{n-2},g_{n-1}^{-1})\\
&=&\gamma_n(g_0,\dots,g_{n-2},g_{n-1}^{-1}g_n)-1,
\end{array}$$ 
so $\gamma_n(h_0,\dots,h_{n-1})-1\in J(R)$ for every $(h_0,\dots,h_{n-1})\models p^{(n)}|\emptyset$. 
But, by Lemma \ref{commutators} and $\omega$-categoricity, there is $l<\omega$ such that every element of $A$ is a product of $l$ elements of the form $\gamma_n(h_0,\dots,h_{n-1})$, where $(h_0,\dots,h_{n-1})\models p^{(n)}|\emptyset$. Also, every element of $R$
 is a sum of a fixed number of elements of the set $A\cup\{-h:h\in A\}\cup\{0\}$. So, we conclude that $J(R)$ is of finite index in $R$.

The rest of the proof follows the lines of the final part of the proof of \cite[Corollary 3.17]{12}. Namely, choose representatives $r_1,\dots,r_m$ of all cosets of $J(R)$ in $R$. By $(*)$, $V_0$ is infinite, and so $R$ and $J(R)$ are infinite as well. Let $k\geq 2$ be the least number such that $J(R)^k=\{0\}$ (such a number exists by Fact \ref{nilp}). Take any non-trivial $i\in J(R)^{k-1}$. Then $i(v)\neq 0$ and $Ai(v)\subseteq Ri(v)=\{r_1i(v),\dots,r_mi(v)\}$. Thus, $Lin_F(Ai(v))$ is a non-trivial, finite-dimensional (over $F$) $FA$-submodule of $\overline{V}$. This is a contradiction with $(*)$, which completes the proof of Claim 1. \hfill $\square$\\

To prove the lemma, we argue by induction on the nilpotency class of $H$. The induction starts by Claim 1. 
Suppose that $H$ is a nilpotent group of class $n+1$ which satisfies the assumptions of the lemma and that the lemma is true for groups of smaller nilpotency class.
Put $C=C_{\Gamma_{n+1}(H)}(V)$. By Claim 1, $C$ is of finite index in $\Gamma_{n+1}(H)$.
Hence, for any $h/C\in \Gamma_{n}(H)/C$ we have that $C_H(h/C)$ is of finite index in $H$, so, by the connectedness of $H$, it is equal to $H$. 
Thus, $H/C$ is nilpotent of class at most $n$. So, by the inductive hypothesis, $C_{H/C}(V)=H/C$. Hence $C_H{V}=H$, as required. \hfill $\blacksquare$\\

Having Fact \ref{solvable} and Lemma \ref{main lemma} in hand, the proof from page 490 of \cite{6} goes through under the assumption of Theorem \ref{Thm1} after a slight modification (which is necessary, because the generic stability of a group is not inherited by definable subgroups).\\

\noindent
{\em Proof of Theorem \ref{Thm1}.} 
Let $G$ be a generically stable group $\emptyset$-definable in a monster model of an $\omega$-categorical theory. By fsg, $G^{00}$ exists, and by $\omega$-categoricity, it is $\emptyset$-definable. Hence, $G^{00}$ has finite index in $G$, and we can assume that $G=G^{00}$. By Fact \ref{solvable}, $G$ is solvable-by-finite, so it has a definable, solvable subgroup of finite index (by $\omega$-categoricity, the group generated by all normal, solvable subgroups of finite index does the job). Thus, $G$ is solvable.

We will argue by induction on the maximal possible length of a chain of distinct, characteristic (in the group-theoretic sense) subgroups of $G$. Suppose $\{e\}=G_0<G_1<\dots<G_t=G$ is a chain of maximal length of distinct, characteristic subgroups of $G$ and that the theorem holds for absolutely connected, generically stable groups with smaller maximal length of such a chain. Notice that all groups $G_i$ are invariant and so $\emptyset$-definable. They are also normal in $G$. 

The group $G_1$ is characteristically simple and solvable, so it is abelian. The group $G/G_1$ is absolutely connected and generically stable, hence, by the induction hypothesis, it is nilpotent-by-finite and so nilpotent (by absolute connectedness). Put $N:=G_1$

Since any nilpotent, $\omega$-categorical group is a direct product of its Sylow subgroups, we may write 
\begin{equation}\tag{$\dag$}
G/N=P_1\times\dots\times P_n, 
\end{equation}
where each $P_i$ is a Sylow $p_i$-subgroup of $G/N$. By $\omega$-categoricity, $G/N$ has bounded exponent, so every $P_i$ is definable. 
Hence, using $(\dag)$, one can easily check that every $P_i$ is a nilpotent, absolutely connected, generically stable group. Thus, we can apply Lemma \ref{main lemma} to definable actions of these groups on the appropriate abelian groups. Having this, the rest of the proof from \cite{6} goes through in our context, which we sketch below.

Let $Q_i$ be the preimage of $P_i$, $i=1,\dots,n$, under the quotient map $G \to G/N$.
Applying Lemma \ref{main lemma} to the action of $P_1$ on the subgroup of $N$ consisting of the elements whose order is co-prime to $p_1$, we get that the elements of $Q_1$ of co-prime orders commute, and so $Q_1$ is locally nilpotent. Since $Q_1$ is also solvable, $\omega$-categoricity and \cite[4(8)]{12} imply that $Q_1$ is nilpotent. Since the group $Q_2Q_1/Q_1$ is definably isomorphic with $P_2/P_2 \cap P_1$, we see that it is a nilpotent, absolutely connected and generically stable $p_2$-group. Thus, applying Lemma \ref{main lemma} to the action of $Q_2Q_1/Q_1$ on successive quotients of the lower central series of the normal Sylow $r$-subgroup of the nilpotent group $Q_1$ (for $r \ne p_2$), we conclude that the elements of co-prime orders in $Q_2Q_1$ commute. As before, this implies that $Q_2Q_1$ is nilpotent. Continuing in this way, we get that $G=Q_n\dots Q_1$ is nilpotent. \hfill $\blacksquare$\\


\noindent
{\bf Address:}\\
Instytut Matematyczny, Uniwersytet Wroc\l awski,\\
pl. Grunwaldzki 2/4, 50-384 Wroc\l aw, Poland.\\[3mm]
{\bf E-mail addresses:}\\
Jan Dobrowolski: dobrowol@math.uni.wroc.pl \\
Krzysztof Krupi\'nski: kkrup@math.uni.wroc.pl


\begin{thebibliography}{99}

\bibitem{13} J. Baldwin, B. Rose, 
{\em $\aleph_0$-categoricity and stability of rings, Journal of Algebra 45, 1-16, 1977.}

\bibitem{1} W. Baur, G. Cherlin, A. Macintyre,
{\em Totally categorical groups and rings}, Journal of Algebra
57, 407-440, 1979. 

\bibitem{15} G. Cherlin, 
{\em On $\aleph_0$-categorical nilrings. II}, Journal of Symbolic Logic 45, 291-301, 1980.

\bibitem{11} J. Dobrowolski, K. Krupi\'nski,
{\em On $\omega$-categorical, generically stable groups}, Journal of Symbolic Logic, accepted.

\bibitem{3} C. Ealy, K. Krupi\'nski, A. Pillay,
{\em Superrosy dependent groups having finitely satisfiable generics}, Annals of Pure and Applied Logic 151,
1-21, 2008.

\bibitem{EW}
D. Evans, F. Wagner, {\em Supersimple $\omega$-categorical groups and theories}, Journal of Symbolic Logic 65, 767-776, 2000.

\bibitem{14}
I. N. Herstein, 
{\em Non-commutative Rings}, Carus Mathematical Monographs, No. 15, Mathematical Association of America, 1968.

\bibitem{7} E. Hrushovski, A. Pillay, 
{\em On NIP and invariant measures}, Journal of the European Mathematical Society, accepted.

\bibitem{4} E. Hrushovski, A. Pillay, Y. Peterzil,
{\em Groups, measures, and the NIP}, Journal of the American Mathematical Society 21, 563-595, 2008. 

\bibitem{5} K. Krupi\'nski,
{\em On $\omega$-categorical groups and rings with NIP}, Proceedings of the American Mathematical Society, accepted.

\bibitem{12} K. Krupi\'nski,
{\em On relationships between algebraic properties of groups and rings in some model-theoretic contexts}, Journal of Symbolic Logic 76, 1403-1417, 2011.

\bibitem{KrWa}
K. Krupi\'nski, F. Wagner, {\em Small profinite groups and rings}, Journal of Algebra 306, 494-506, 2006.

\bibitem{6} H. D. Macpherson, {\em Absolutely ubiquitous structures and $\aleph_{0}$-categorical groups}, Quart. J.
Math. Oxford (2) 39, 483-500, 1988. 

\bibitem{8} A. Pillay, P. Tanovi\'c, {\em Generic stability, regularity, and quasi-minimality},
CRM Proceedings and Lecture Notes 53, 189-211, 2011.

\bibitem{9} B. Poizat,
{\em Stable groups}, American Mathematical Society, Providence, 2001.


\bibitem{Wa}
F. Wagner, {\em Stable groups}, London Math. Soc. Lecture Notes Series 240, Cambridge University Press, UK, 1997.

\bibitem{10} J. Wilson,
{\em The algebraic structure of $\omega$-categorical groups}, in: Groups-St. Andrews, Ed. C.
M. Campbell, E. F. Robertson, London Math. Soc. Lecture Notes 71, Cambridge, 345-358,
1981. 


\end{thebibliography}
\end{document}